\documentclass{amsart}

\usepackage{booktabs}
\usepackage{multirow}
\usepackage{array}
\usepackage{amsthm,amssymb,latexsym,amsmath}
\usepackage{mathrsfs}
\usepackage{graphicx}
\usepackage{color}
\usepackage[brazil]{babel}
\usepackage[latin1]{inputenc}
\usepackage{graphicx}%
\usepackage{amsmath,amssymb,amsfonts}%
\usepackage{amsthm}%
\usepackage{physics}
\usepackage{mathrsfs}%
\usepackage[title]{appendix}%
\usepackage{xcolor}%
\usepackage{textcomp}%
\usepackage{manyfoot}%
\usepackage{booktabs}%
\usepackage{algorithm}%
\usepackage{algorithmicx}%
\usepackage{algpseudocode}%
\usepackage{listings}%
\usepackage{amsthm,amssymb,latexsym,amsmath}
\usepackage{mathrsfs}
\usepackage{graphicx}
\usepackage{color}
\usepackage[all]{xy}
\usepackage{graphicx}%
\usepackage{multirow}%
\usepackage{amsmath,amssymb,amsfonts}%
\usepackage{amsthm}%
\usepackage{mathrsfs}%
\usepackage[title]{appendix}%
\usepackage{xcolor}%
\usepackage{textcomp}%
\usepackage{manyfoot}%
\usepackage{booktabs}%
\usepackage{algorithm}%
\usepackage{algorithmicx}%
\usepackage{algpseudocode}%
\usepackage{listings}%

\usepackage{hyperref}

\newcommand{\CC}{{\mathbb C}}

\renewcommand{\dim}{\mathrm{dim}}

\newcommand{\OO}{\mathcal O}

\newcommand{\fol}{\mathcal{F}}

\newcommand{\NN}{{\mathbb{N}}}

\newcommand{\PP}{\mathbb{P}}

\newcommand{\reg}{\mathrm{reg}\,}

\newcommand{\Sing}{{\rm Sing}}

\newtheorem{lema}{Lemma}[section]
\newtheorem{cor}[lema]{Corollary}

\newtheorem*{teo1234}{Theorem}
\newtheorem*{teo1}{Theorem 1}
\newtheorem*{teo2}{Theorem 2}
\newtheorem*{teo3}{Theorem 3}
\newtheorem*{teo4}{Theorem 4}
\newtheorem*{teo5}{Theorem 5}

\theoremstyle{definition}
\newtheorem{remark}[lema]{Remark}

\newtheorem{exe}[lema]{Example}

\begin{document}

\title{Bounds and formulas for residues of singular holomorphic foliations and applications}

\dedicatory{\it Dedicated to Catarina Vitória}

\begin{abstract}
We consider one dimensional holomorphic foliations with isolated singularities that leave invariant a local complete intersection. We establish explicit formulas for the total GSV index of such foliations and obtain bounds for this index. As applications, we derive several consequences related to Poincaré's problem for foliations on projective spaces.
\end{abstract}

\author{Diogo da Silva Machado}

\address{\noindent  Diogo da Silva Machado\\
Departamento de Matem\'atica \\
Universidade Federal de Vi\c cosa\\
Avenida Peter Henry Rolfs, s/n - Campus Universitário \\
36570-900 Vi\c cosa- MG, Brazil} \email{diogo.machado@ufv.br}




\maketitle



\section{Introduction}

The index of a vector field at a singular point is an elementary and fundamental notion in Geometry and Topology, having been extensively developed and studied from various perspectives. Its properties underlie crucial results, most notably the classical Poincaré-Hopf theorem, which states that the sum of the indices of a vector field with isolated singularities on a closed, smooth, oriented manifold is a topological invariant, coinciding with the Euler characteristic of the manifold.

However, when considering spaces that are not manifolds --- such as singular varieties --- it is natural to seek a generalization of this index. In response, X. Gómez-Mont, J. Seade, and A. Verjovsky \cite{GSV} introduced the so-called GSV index. This concept extends the classical Poincaré - Hopf index while crucially preserving its stability under perturbations. Initially defined for vector fields on hypersurfaces with isolated singularities \cite{GSV}, the GSV index was later extended to complete intersections in \cite{SeaSuw1}, and further generalizations have since been obtained in various contexts by several authors \cite{a03}, \cite{a09}, \cite{cm3}, \cite{a08}, \cite{a07}.

By employing a suitable adaptation of Chern-Weil theory to \v{C}ech-de Rham cohomology, the Baum-Bott residue theory \cite{BB1} has been extended to the setting of foliations restricted to their invariant subvarieties. In this framework, the GSV index is defined as a residue of the foliation relative to a singular invariant subvariety, and is obtained as the localization of a certain virtual bundle associated with both the foliation and the invariant subvariety (see, for example, \cite{Bru}, \cite{MD3}, \cite{LS01}, \cite{DSM}, \cite{DSM2} and \cite{SeaSuw2}). In general, the computation of these residues is a challenging task, and in some cases, alternative approaches have been employed to compute them. For example, Y. Nakamura (Kinki University, Japan) used computer software to compute the GSV index of certain holomorphic vector fields (see \cite[Remark 5.7.1.]{BarSaeSuw}). 

We refer to Section \ref{sec111} for more details on the GSV index in the context of one-dimensional holomorphic foliations. 

Let $\fol$ be a one-dimensional holomorphic foliation, with isolated singularities, on a compact complex manifold $X$, and let $V\subset X$ be a variety, with isolated singularities, such that $\fol$ leaves $V$ invariant.  For each singular point $p \in S(\fol, V) := (Sing(\fol) \cap V) \cup Sing(V)$, we shall denote by $GSV_p(\fol,V)$ the GSV index of $\fol$ at $p$ along $V$, and we shall write $GSV(\fol,V)$ for the sum $\sum_{p\in S(\fol,V)} GSV_p(\fol,V)$.

In the present work, we establish a formula for computing $GSV(\fol,V)$, which
 will also be called the {\it total GSV index of $\fol$}. We consider the case where $V$ is a local complete intersection defined as the zero set of a section of a vector bundle (this is the case, for example, of hypersurfaces, global complete intersections and complete intersections in projective space). More precisely, we prove the following formula:

\begin{teo1}\label{ttt1}
Let $X$ be an $m$-dimensional complex compact manifold and $\fol$ a one-dimensional holomorphic foliation on $X$, with isolated singularities, leaving invariant a local complete intersection $V$, of codimension $r$, with isolated singularities, defined by a section of a holomorphic
vector bundle $N$ on X. We have,

{\footnotesize {
\begin{eqnarray} \label{GSVF}\nonumber
GSV(\fol,V) \hspace{-0.07cm} =  \hspace{-0.15cm} \int_X \hspace{-0.1cm} c_r(N)\hspace{-0.12cm}\left(\sum_{t=0}^{m-r}\hspace{-0.14cm}\left(\hspace{-0.1cm} c_t(TX) \hspace{-0.05cm} + \hspace{-0.09cm} \sum_{j=1}^t\hspace{-0.05cm}\sum_{i=1}^j\hspace{-0.05cm}\sum_{\mid L_i\mid=j}\hspace{-0.19cm}c_{t-j}(TX)(-1)^i\hspace{-0.05cm}c_{L_i}\hspace{-0.14cm}\left(N\right)\hspace{-0.04cm} \right)\hspace{-0.1cm}c_1(T\fol^{\ast})^{m-r-t}\right) 
\end{eqnarray}
}}

\medskip

\noindent where $T\fol^{\ast}$  denotes the dual tangent bundle of the $\fol$, $L_i = (l_1,\ldots,l_i)\in (\NN^{\ast})^i$, with  $l_1+\ldots +l_i = j$, is an $i$-dimensional multi-index and $\displaystyle c_{L_i}\left(N\right) := c_{l_1}\left(N\right)\cdots c_{l_i}\left(N\right)$.
\end{teo1}

Bounds for the GSV index are established in Theorem 2, in the case where $\fol$ has non-degenerate isolated singularities. Recall that $p\in Sing(\fol)$ is a non-degenerate singularity when the differential $Dv(p)$ of $v$ at $p$ is non-singular, where $v$ is a holomorphic vector field that defines $\fol$ locally at $p$. 

\begin{teo2}\label{ttt}
Let $\fol$ be a one-dimensional holomorphic foliation, with non-degenerate isolated singularities, over an $m$-dimensional complex manifold $X$, and $V\subset X$ a local complete intersection of codimension $r$, with isolated singularities, invariant by $\fol$.  For all $p\in Sing(\fol) \cap Sing(V)$ one has

\medskip

\noindent (i) If $\dim(V)$ is even then

{\small{ 
$$
\alpha + \tau_p(V) \,\,\, \leq \,\,\, GSV_p(\fol, V) \,\,\, \leq \,\,\hspace{-0.09cm} \beta + \tau_p(V)
$$
}}

\bigskip

\noindent (ii) If $\dim(V)$ is odd then

{\small{ 
$$
\beta - \tau_p(V) \,\,\, \leq \,\,\, GSV_p(\fol, V) \,\,\, \leq \,\, \hspace{-0.09cm} \alpha - \tau_p(V)
$$
}}

\noindent where $\alpha$ and $\beta$ are integer constants given by
$$
\alpha = \sum_{j=0}^{m-r-1} (-1)^j {r-1+j \choose j},\,\,\,\,\,\,\, \beta = \alpha + (-1)^{m-r-1}{m-2 \choose m-r-1}
$$ 
\noindent and $\tau_p(V)$ denotes the invariant of Greuel of the isolated complete intersection singularity germ $(V,p)$.
\end{teo2}

\begin{remark} \label{rr1} Recall that if $f = (f_1,\ldots, f_r) = 0$ is a local reduced equation of $(V,p)$, in local coordinates $(x_1,\ldots, x_m)$, the invariant of Greuel $\tau_p(V)$ can be calculated by the following formula (see \cite{ggel})

$$
\tau_p(V) = \displaystyle dim_{\CC}\frac{\OO_{\CC^m,p}}{\left\langle f,\displaystyle\frac{\partial(f_1,\ldots,f_r)}{\partial(x_{j_1},\ldots,x_{j_r})} (1\leq j_1<\ldots<j_r \leq m)\right\rangle},
$$
\noindent where $\OO_{\CC^m,p}$ denotes the ring of germs of holomorphic functions at $(\CC^m,p)$. When $(V,p)$ is a hypersurface germ, the invariant of Greuel $\tau_p(V)$ is also referred to as the {\it Tjurina number} (see \cite{ggel}).
\end{remark}

\hyphenation{di-ffe-ren-tial}
\hyphenation{inva-riant}
\hyphenation{sin-gu-la-ri-ties}

As an application of Theorem 1, we obtain some results related to Poincaré problem. We recall that the Poincaré problem is the question of bounding the degree of an algebraic variety invariant by a holomorphic foliation $\fol$ of complex projective space $\PP^m$, in terms of the degree of the foliation. This question was first considered by Poincaré in \cite{Poin}, in the case of plane curves invariant by foliations of $\PP^2$ with a rational first integral. More recently, D. Cerveau and A. Lins Neto \cite{Alc}, considering foliations on the projective space $\PP^2$, proved that if an $\fol$-invariant curve $C$ has at most nodal singularities then 
\begin{eqnarray}\label{scv1}
k \leq  d + 2,
\end{eqnarray}

\noindent where $k$ and $d$ denotes, respectively, the degrees of $C$ and $\fol$. M. Carnicer \cite{Can}  has obtained the same inequality when the singularities of $\fol$ over $C$ are all non-dicritical. M. Brunella \cite{Bru} shows that the negativity of the GSV-indices is an obstruction to the inequality (\ref{scv1}). He recovered Carnicer's result by observing that these indices are always non-negative in the non-dicritical case. 

In \cite{DSM}, we show that the GSV index remains an obstruction to inequality (\ref{scv1}) when $\fol$ is a foliation on a projective space of even dimension, leaving invariant a hypersurface with isolated singularities (in odd dimension, we also prove that the sum of the GSV indices is always non-negative; nevertheless, this is not sufficient to guarantee that inequality (\ref{scv1}) holds). 

In the present work, we aim to show that, in the case where $\fol$ leaves invariant a complete intersection curve, this result holds for projective spaces of arbitrary dimension. More precisely:

\begin{teo3}\label{tectec}
Let $\fol$ be a holomorphic foliation of dimension one in $\PP^m$ and degree $d$, with isolated singularities, leaving invariant a complete intersection curve $C \subset \PP^m$ of multi-degre $(k_1,\ldots, k_{m-1})$. Then inequality
$$
k_1 + \ldots +k_{m-1} \leq d + m.
$$
\noindent holds if and only if $\displaystyle GSV(\fol,V) \geq 0$.
\end{teo3}

\begin{remark} {\bf (a)} A complete intersection curve $C$ on $\PP^m$ is defined by a section of vector bundle $\OO_{\PP^m}(1)^{\otimes k_1}\oplus \ldots \oplus \OO_{\PP^m}(1)^{\otimes k_{m-1}}$, where $\OO_{\PP^m}(1)$ denotes the hyperplane bundle. We recall that the multi-degree of $C$ is defined by $(k_1,\ldots, k_{m-1})$.\\\\
\noindent {\bf (b)} In \cite{CL}, V. Cavalier and D. Lehmann, using a bound for total GSV index, proved a type of Poincaré inequality for one-dimensional foliations on the projective space.
\end{remark}

An immediate consequence of the Theorem 3 is the following.

\begin{cor}[M. G. Soares, \cite{SoaAnn}] 
If $\fol$ is a foliation of dimension one in $\PP^m$ and degree $d$, with isolated singularities, leaving invariant a smooth complete intersection curve $C \subset \PP^m$ of multi-degre $(k_1,\ldots, k_{m-1})$, then 
$$
k_1 + \ldots +k_{m-1} \leq d + m.
$$
\end{cor}

In the study of the Poincaré problem, by considering the degree and other characters associated with an invariant curve, M. G. Soares \cite{Soa} also proved the following Poincaré inequality for holomorphic foliations on the two-dimensional projective space:

\medskip

\begin{teo1234}[M. G. Soares, \cite{Soa}]
Let $\fol$ be a holomorphic  foliation of dimension one in $\PP^2$, with isolated singularities and degree $d$, leaving invariant a curve $C$ of degree $k$. One has,
\begin{eqnarray}\label{9}
k(k - 2) \hspace*{0.25cm} - \sum_{p\in Sing(C)}(\mu_p(C)- 1) \hspace*{0.25cm}\leq \hspace*{0.25cm}\displaystyle dk. 
\end{eqnarray}
\noindent where $\mu_p(C)$ denotes the Milnor number of $C$ at $p$.
\end{teo1234}

\noindent   In his proof, M. G. Soares also considered $Sing(C)\subset Sing(\fol)$. We establish a generalization of the formula (\ref{9}) for the case of arbitrary dimension $m$:

\begin{teo4}\label{ecd}
Let $\fol$ be a holomorphic foliation of dimension one on $\PP^m$, with isolated singularities and degree $deg(\fol) = d$, leaving invariant a complete intersection curve $C$ of multi-degree $(k_1,\ldots,k_{m-1})$, with $Sing(C)\subset Sing(\fol)$. One has,
$$
(k_1 \ldots k_{m-1}) \left(k_1 +\ldots + k_{m-1} - m \right)  \hspace{0.0cm} - \sum_{p\in Sing(C)}\hspace{-0.35cm}(\mu_p(C) -1) \leq d(k_1 \ldots k_{m-1}),
$$
\noindent where $\mu_p(C)$ denotes the Milnor number of $C$ at $p$.
\end{teo4}

\begin{remark} In \cite{DSM} we obtained a characterization of Soares's inequality (\ref{9}) in terms of the GSV index, considered hypersurfaces invariant by holomorphic foliations on projective space of arbitrary dimension. In \cite{DSM2}, generalization of (\ref{9}) was established, addressing the case of holomorphic foliations on projective spaces of even dimension.
\end{remark}

In Section \ref{prel}, we fix our notation and recall some basic facts, including the definition of the GSV index in the context of one-dimensional holomorphic foliations. Theorem 1 will be proven in Section \ref{secc992}, where other immediate consequences can also be found (Corollary \ref{coroo} and Corollary \ref{GSformula}). In Section \ref{secc03}, we devote ourselves to establishing bounds for the GSV index, culminating in the proof of Theorem 2. 

In Section \ref{secc04} we turn to projective foliations with the aim of obtaining applications in the context of the Poincaré Problem. Our goal is to present the proofs of Theorem 3 and Theorem 4. Moreover, we shall present a few remarks  on Schwartz index in the context of holomorphic foliations. 

\section{Preliminaries}\label{prel}
\subsection{Singular one-dimensional holomorphic foliations} Let $X$ be a connected complex manifold. A one-dimensional holomorphic foliation $\fol$ on $X$ is given by a set of data $\{(U_{\alpha}, v_a, g_{\alpha\beta})\}$, where $\{U_{\alpha}\}$ is an open covering of $X$, $v_\alpha \in \mathrm{H}^0(U_{\alpha},TX|_{U_{\alpha}})$ is a holomorphic vector field on $U_{\alpha}$, for each $\alpha$, and $g_{\alpha\beta} \in \mathcal{O}_X^*(U_\alpha\cap U_\beta)$ is a non-vanishing holomorphic function on $U_{\alpha} \cap U_{\beta}$, for each pair $(\alpha,\beta)$, such that $v_\alpha = g_{\alpha\beta}v_\beta$ in $U_\alpha\cap U_\beta$ and $g_{\alpha\beta}g_{\beta\gamma} = g_{\alpha\gamma}$ in $U_\alpha\cap U_\beta\cap U_\gamma$.

We denote by $K_{\fol}$ the line bundle defined by the cocycle $\{g_{\alpha\beta}\}\in \mathrm{H}^1(X, \mathcal{O}^*)$. Thus, a one-dimensional holomorphic  foliation $\fol$ on $X$ induces  a global holomorphic section $s_{\fol}\in \mathrm{H}^0(X,T_X\otimes K_{\fol})$. The line bundle $T_{\fol}:= (K_{\fol})^* \hookrightarrow T_X$ is called the tangent  bundle of $\fol$. The singular set of $\fol$ 
is $Sing(\fol)=\{s_{\fol}=0\}$. We will assume that $codim (Sing(\fol))\geq 2$.

Given a variety $V \subset X$, we say that the foliation $\fol$ \emph{leaves} $V$ \emph{invariant} (or that $V$ is {\it invariant by $\fol$} ) if it satisfies the following condition: for all $x\in V - Sing(V)$, the vector $v_{\alpha}(x)$ belongs to $T_xV$, with $x\in U_{\alpha}$. 

A foliation on a complex projective space $\mathbb{P}^n$  is  called a {\it projective}. Let $\fol $ be  a   projective foliation   with tangent bundle $T_{\fol}=\mathcal{O}_{\mathbb{P}^n}(r)$. The integer $d:=1-r$ is called the degree of $\fol$.

\subsection{The GSV-index} \label{sec111} X. Gomez-Mont, J. Seade and A. Verjovsky \cite{GSV} introduced the GSV-index for a holomorphic vector field over an analytic hypersurface with isolated singularities, on a complex manifold, generalizing the (classical) Poincar\'e-Hopf  index. The concept of GSV-index was extended to continuous vector fields on more general contexts. For example, J. Seade and T. Suwa in \cite{SeaSuw1}, defined the GSV-index for vector fields on isolated complete intersection singularity germs. J.-P. Brasselet, J. Seade and T. Suwa in \cite{a03}, extended the notion to vector fields defined on certain types of analytic subvariety with non-isolated singularities. 

\hyphenation{using}
In \cite{a08}   X. Gomez-Mont  defined the homological index of holomorphic vector fields on an analytic hypersurface with isolated singularities, which coincides with the GSV index. There is also the virtual index, introduced by D. Lehmann, M. Soares and T. Suwa \cite{a07}, defined via Chern-Weil theory, that can be interpreted as the GSV-index. 
M. Brunella \cite{a09} also presents the GSV-index for foliations on complex surfaces  by a different approach.

\hyphenation{sin-gu-la-ri-ties}
\hyphenation{sin-gu-la-ri-ty}
\hyphenation{de-fi-ning}
\hyphenation{fields}

Let us recall the definition of the GSV index (see, for example, \cite[Ch.3, Sec. 3.2]{BarSaeSuw}. Let $X$ be an $m$-dimensional complex manifold and $\fol$ a one-dimensional foliation on $X$, with isolated singularities, leaving invariant a local complete intersection $V$, of codimension $r$, with isolated singularities. Given $p_0 \in Sing(V)$, let $f = (f_1,\ldots, f_r) = 0$ be a local reduced equation of $V$ at $p_0$ and  $v$ a holomorphic vector field defining $\fol$ on a neighborhood $U$ of $p_0$ in $X$. The gradient vector fields $\{\overline{grad}\,(f_1), \ldots, \overline{grad}\,(f_r) \}$ are linearly independent everywhere away from $p_0$, because $p_0$ is an isolated singularity of $V$. 
Since $\fol$ leaves $V$ invariant, it follows that the restriction $v_{\ast}$ of vector field $v$ to the  regular part $V_{reg} \cap U$ defines a continuous vector field on $V\cap U$, possibly singular only in $p_0$.

Let $z=(z_1,\ldots,z_m)$ be a coordinate system on $U$ and consider 
$$
S_{\varepsilon} = \{z= (z_1,\ldots,z_m):\,\,\mid\mid z-p_0\mid\mid \,\, = \,\, \varepsilon\}$$  
\noindent a sphere with $\varepsilon$ sufficiently small so that $K =  V \cap S_{\varepsilon} $ is the link of the singularity of $ V$ at $p_0$ (see, for example, \cite{JM}). It is an $(2m-1)$-dimensional real oriented manifold. The set $\{v_{\ast}(z), \overline{grad}\,(f_1)(z), \ldots, \overline{grad}\,(f_r)(z) \}$ is a $(r+1)$-frame at each point $z\in V - \{p_0\}$, and up to homotopy, it can be assumed to be orthonormal. Thus, we obtain a continuous map
$$
\phi_{v_{\ast}}:=(v_{\ast}(z), \overline{grad}\,(f_1)(z), \ldots, \overline{grad}\,(f_r)(z)): K \longrightarrow W_{r+1}(m+r)
$$  
\noindent where $W_{r+1}(m+r)$ is the Stiefel manifold of complex $(r+1)$-frames in $\CC^{m+1}$. The GSV index of $v$ at $p_0 \in V$, denoted by $GSV_{p_0}(v, V)$, is defined as the degree of map $\phi_{v_{\ast}}$. According to \cite[Ch.6, Corollary 6.3.1]{BarSaeSuw}, $GSV_{p_0}(v, V)$ does not depend on the choice of the vector field $v$, local representative of the foliation $\fol$ at $p_0$. Thus, the {\it GSV index of $\fol$ in $p_0$ along $V$}, denoted by $GSV_{p_0}(\fol, V)$, is defined as $GSV_{p_0}(v, V)$. Moreover, if the local
complete intersection $V$ is defined by a section of a holomorphic vector bundle $N$ on X, then we have the following residue type formula for total GSV index of the foliation (see \cite[Ch.IV, Theorem 7.16]{Suw2}) 
$$
\displaystyle\int_{V} c_{n-r}(TX - N- T_{\fol}) = GSV(\fol,V).
$$

\medskip

\begin{remark} In the above definition of the GSV index, we suppose that $V$ has dimension greater than 1 or it has only one branch. When $dim(V) =1$ and $V$ has several branches the GSV index is defined as the Poincaré-Hopf index of an extension of $v$ to a Milnor fiber (see \cite{Bru1}, \cite{Bru}, \cite{KS}). Also, we observe that if $p_0 \in \Sing(\fol)\cap V_{\reg}$ then $GSV_{p_0}(\fol, V)$ is the (classical) Poincaré-Hopf index of $v_{\ast}$ at $p_0$.
\end{remark}

\medskip

\section{Formulas for GSV index}\label{secc992}

In this section, we present the proof of Theorem 1, which provides an explicit formula for the total GSV index of the foliation. We also provide versions for this formula for the case in where the invariant variety is a complete intersection of hypersurfaces (Corollary \ref{coroo}), including the case of projective foliations (Corollary \ref{GSformula}).

To prove Theorem \ref{ttt1}, we will first prove the following result.
\begin{lema}\label{proteo}
Let $X$ be an $m$-dimensional complex manifold and $V\subset X$ a local complete intersection of codimension $r$, with isolated singularities, defined by a section of a holomorphic vector bundle $N$. Then, for all $t=1,\ldots, m$, the $t$-th Chern class of the virtual bundle $TX-N$ is given by {\small{
\begin{eqnarray}\label{revrev00}
c_t\left(TX - N\right) = c_t(TX) + \sum_{j=1}^t\sum_{i=1}^j\sum_{\mid L\mid=j}(-1)^ic_{t-j}(TX)c_{L_{i}}\left(N\right),
\end{eqnarray}
}}
\noindent where $L_i = (l_1,\ldots,l_i)\in (\NN^{\ast})^i$, with  $l_1+\ldots +l_i = j$, is an $i$-dimensional multi-index and $ c_{L_i}\left(N\right) := c_{l_1}\left(N\right)\cdots c_{l_i}\left(N\right)$               
\end{lema}

\noindent {\bf Proof.:} By definition the total Chern class of the virtual vector bundle $TX-N$ is $c_{\ast}(TX-N) = c_{\ast}(TX).c_{\ast}(N)^{-1}$ (see \cite{BB1}).  Thus, we obtain the relations               
{\small{
\begin{eqnarray}\nonumber
\left\{
\begin{array}{lcc}
c_1(TX-N) &=& c_1(TX) - c_1(N)\\
c_2(TX-N) + c_1(N)c_1(TX-N) &=& c_2(TX) - c_2(N)\\ \vdots & \vdots&\vdots \\
c_{t}(TX-N) + c_1(N)c_{t-1}(TX-N) + \ldots + c_{t-1}(N)c_{1}(TX-N)  &=& c_t(TX) - c_t(N)
\\
\end{array} \right.
\end{eqnarray}
}}

\noindent which can be written in the following matrix identity

$${\scriptsize{
{\left(
\begin{array}{cccccccccccccc}
1 & 0 & 0 & \ldots & 0 & 0 \\
c_1(N) & 1 & 0 & \ldots & 0 & 0 \\
c_2(N) & c_1(N) & 1 & \ldots & 0 & 0 \\
\vdots & \vdots & \vdots & \ddots & \vdots & \vdots\\c_{t-2}(N) & c_{t-3}(N) & c_{t-4}(N) & \ldots & 1 & 0\\
c_{t-1}(N) & c_{t-2}(N) & c_{t-3}(N) & \ldots & c_1(N) & 1
\end{array}
\right)\left(\begin{array}{cccccccccccccc}
c_1(TX-N)\\
c_2(TX-N)\\
c_3(TX-N)\\
\vdots \\
c_{t-1}(TX-N)\\
c_{t}(TX-N)
\end{array}
\right) = \left(\begin{array}{cccccccccccccc}
c_1(TX) - c_1(N)\\
c_2(TX) - c_2(N)\\
c_3(TX) - c_3(N) \\
\vdots \\
c_{t-1}(TX) - c_{t-1}(N)\\
c_{t}(TX) - c_t(N)
\end{array}
\right)}}}$$

\medskip

\noindent or, equivalently,

$${\scriptsize \left(
\begin{array}{cccccccccccccc}
1 &\hspace{0.9cm} 0 \hspace{0.9cm} & 0 & \ldots & 0 & \hspace{0.18cm} 0 \hspace{0.18cm} \\
0 & 1 & 0 & \ldots & 0 & 0 \\
0 & 0 & 1 & \ldots & 0 & 0 \\
\vdots & \vdots & \vdots & \ddots & \vdots & \vdots\\0 & 0 & 0 & \ldots & 1 & 0\\
0 & 0 & 0 & \ldots & 0 & 1
\end{array}
\right)\left(\begin{array}{cccccccccccccc}
c_1(TX-N)\\
c_2(TX-N)\\
c_3(TX-N)\\
\vdots \\
c_{t-1}(TX-N)\\
c_{t}(TX-N)
\end{array}
\right) = \left(\begin{array}{cccccccccccccc}
\delta_1(TX,N)\\
\delta_2(TX,N)\\
\delta_3(TX,N) \\
\vdots \\
\delta_{t-1}(TX,N)\\
\delta_t(TX,N)
\end{array}
\right)}$$

\medskip

\noindent where the $\delta_j(TX,N)'$s are defined by the recursion 
\begin{eqnarray}\nonumber\label{recu}
\left\{
\begin{array}{lcl}
\delta_1(TX,N) &=& c_1(TX) - c_1(N)\\ 
\delta_j(TX,N) &=& c_{j}(TX) - c_{j}(N) - \displaystyle\sum_{i=1}^{j-1}c_{j-i}(N)\delta_i(TX,N),\,\,\, \mbox{for $j=2,3,\ldots,k$}. 
\end{array} \right.
\end{eqnarray}

\noindent The verification of this equivalence consists of applying a sequence of $(t-1)!$ elementary row operations, divided into $t-1$ steps. The operations corresponding to the $j$-th step ($j = 1, 2, \ldots, t-1$) are defined in each row $R_{\ast}$ by:$$R_{j+1} \leftarrow R_{j+1} - c_{1}(N)R_{j}, \quad
R_{j+2} \leftarrow R_{j+2} - c_{2}(N)R_{j}, \quad \ldots, \quad
R_{t} \leftarrow R_{t} - c_{t-j}(N)R_{j}.
$$

\noindent Formula (\ref{revrev00}) is then obtained by successively applying recursion in the relation $c_{t}(TX-N) = \delta_t(TX,N)$.

\noindent $\square$

\noindent {\bf Proof of Theorem 1:} By hypothesis the $\fol$-invariant variety $V$ is
a local complete intersection defined as the zero set of a holomorphic section $\upsilon \in H^0(X,N)$, which is regular (see, for example, \cite[Ch.II, Sec. 13]{Suw2}). Hence, the $r$-th Chern class of $N$ is Poincaré dual to the fundamental class of $V$ and, consequently,
 
\begin{eqnarray} \nonumber
&&\hspace{-0.7cm}\int_X c_r(N)\wedge\sum_{t=0}^{m-r}\left( c_t(TX) + \sum_{j=1}^t\sum_{i=1}^j\sum_{\mid L_i\mid=j}(-1)^ic_{t-j}(TX)c_{L_i}\left(N\right)\right)c_1(T\fol^{\ast})^{m-r-t} = \\\nonumber && \hspace{-0.7cm}
 = \int_{V} \sum_{t=0}^{m-r}\left( c_t(TX) + \sum_{j=1}^t\sum_{i=1}^j\sum_{\mid L_i\mid=j}(-1)^ic_{t-j}(TX)c_{L_i}\left(N\right)\right)c_1(T\fol^{\ast})^{m-r-t}.
\end{eqnarray}

\noindent By using the Lemma \ref{proteo} we compute the last sum as follows:

\begin{eqnarray} \nonumber
&& \int_{V} \sum_{t=0}^{n-r}\left( c_t(TX) + \sum_{j=1}^t\sum_{i=1}^j\sum_{\mid L_i\mid=j}(-1)^ic_{t-j}(TX)c_{L_i}\left(N\right)\right)c_1(T\fol^{\ast})^{m-k-t} = \\\nonumber &&= \int_{V}  \sum_{t=0}^{m-r}c_{t}\left(TX - N\right) c_1(T\fol^{\ast})^{m-r-t}= \\\nonumber &&= \int_{V}  c_{m-r}\left(TX - N - T\fol \right) = \\\nonumber &&= GSV(\fol,V)
\end{eqnarray}

\noindent where in the last step we are using \cite[Ch.IV, Theorem 7.16]{Suw2}.

\noindent $\square$

\begin{cor}\label{coroo}
Let $X$ be an $m$-dimensional complex compact manifold and $\fol$ a one-dimensional foliation on $X$, with isolated singularities, leaving invariant the complete intersection $V = \displaystyle \cap_{\lambda =1}^rD_{\lambda}$ of $r$ hypersurfaces $D_{\lambda}$. Suppose $V$ has isolated singularities. One has, 
{\footnotesize {
\begin{eqnarray} \label{GSVF}\nonumber
GSV(\fol,V) \hspace{-0.07cm} =  \hspace{-0.15cm} \int_X \hspace{-0.1cm} \prod_{\lambda = 1}^r\hspace{-0.1cm}c_1([D_{\lambda}])\hspace{-0.16cm}\sum_{t=0}^{m-r}\hspace{-0.14cm}\left(\hspace{-0.1cm} c_t(TX) \hspace{-0.05cm} + \hspace{-0.09cm} \sum_{j=1}^t\hspace{-0.05cm}\sum_{i=1}^j\hspace{-0.05cm}\sum_{\mid L_i\mid=j}\hspace{-0.19cm}c_{t-j}(TX)(-1)^i\hspace{-0.05cm}c_{L_i}\hspace{-0.14cm}\left(\bigoplus^r_{\lambda=1}[D_{\lambda}]\right)\hspace{-0.14cm} \right)\hspace{-0.1cm}c_1(T\fol^{\ast})^{m-r-t} 
\end{eqnarray}
}}
\noindent where  $[D_{\lambda}]$ denotes the line
bundle of divisor $D_{\lambda}$,  $L_i = (l_1,\ldots,l_i)\in (\NN^{\ast})^i$, with  $l_1+\ldots +l_i = j$, is an $i$-dimensional multi-index and\\ 
\noindent $\displaystyle c_{L_i}\left(\bigoplus^r_{\lambda=1}[D_{\lambda}]\right) := c_{l_1}\left(\bigoplus^r_{\lambda=1}[D_{\lambda}]\right)\cdots c_{l_i}\left(\bigoplus^r_{\lambda=1}[D_{\lambda}]\right)$.
\end{cor}

\noindent {\bf Proof.:} Covering $X$ by open sets $U_{\alpha}$ so that each hypersurface $D_{\lambda}$ is defined by equation $f_{\lambda ,\alpha} = 0$, in $U_{\alpha}$, we have that each holomorphic bundle $[D_{\lambda}]$ admits a natural holomorphic section determined by the collection $\{f_{\lambda, \alpha}\}$. The collection $\{\bigoplus^r_{\lambda=1}f_{\lambda \alpha}\}$ determines a holomorphic section $\upsilon$ of $\bigoplus^r_{\lambda=1}[D_{\lambda}]$ whose zero set $Z(\upsilon)$ is exactly $V = \cap^r_{\lambda=1}D_{\lambda}$. Since $c_r(\bigoplus^r_{\lambda=1}[D_{\lambda}])=\prod_{\lambda = 1}^rc_1([D_{\lambda}])$  it follows from Theorem 1 that

{\small {
\begin{eqnarray} \nonumber
&& \hspace{-0.9cm}GSV(\fol,V) = \\\nonumber && \hspace{-0.7cm}=\int_X \hspace{-0.1cm} \prod_{\lambda = 1}^r\hspace{-0.1cm}c_1([D_{\lambda}])\hspace{-0.1cm}\sum_{t=0}^{n-r}\hspace{-0.1cm}\left(\hspace{-0.1cm} c_t(TX) \hspace{-0.1cm} + \hspace{-0.1cm} \sum_{j=1}^t\sum_{i=1}^j\sum_{\mid L_{i}\mid=j}\hspace{-0.1cm}(-1)^ic_{t-j}(TX)c_{L_i}\hspace{-0.1cm}\left(\bigoplus^r_{\lambda=1}[D_{\lambda}]\hspace{-0.1cm}\right)\hspace{-0.1cm}\right)c_1(T\fol^{\ast})^{m-r-t} \hspace{-0.1cm} \end{eqnarray}
}}

\noindent $\square$

In the case where $\fol$ is a one-dimensional holomorphic foliation on the projective space $\PP^m$, we proved Corollary \ref{GSformula} which establishes formula for $GSV(\fol,V)$ in terms of the degree of the foliation and the degree of the invariant  complete interserction $V$. More precisely, one has:
  
\begin{cor}\label{GSformula}
Let $\fol$ be a foliation of dimension one on $\PP^n$, with isolated singularities and degree $\deg(\fol) = d$, leaving invariant the complete intersection $V =  \bigcap_{\lambda =1}^rD_{\lambda}$ of $r$ hypersurfaces $D_{\lambda}$, each of degree $\deg(D_{\lambda})$. Suppose that $V$ has isolated singularities. One has, 

{\footnotesize {
\begin{eqnarray} \label{GSVF}\nonumber
\hspace{-0.15cm} GSV(\fol,V) \hspace{-0.07cm} = \hspace{-0.13cm} \prod_{\lambda = 1}^r \hspace{-0.1cm} deg(D_{\lambda})\hspace{-0.15cm}\left[\sum_{t=0}^{m-r}\hspace{-0.1cm}\left( \hspace{-0.15cm} {m+1 \choose t} \hspace{-0.1cm} + \hspace{-0.12cm} \sum_{j=1}^t\sum_{i=1}^j\sum_{\mid L_i\mid=j}\hspace{-0.19cm}(-1)^i\hspace{-0.05cm}{m+1 \choose t-j}deg_{L_i}\hspace{-0.09cm}\left[D_1,\ldots, D_r\right]\hspace{-0.1cm}\right)\hspace{-0.15cm}(d-1)^{m-r-t} \hspace{-0.08cm} \right] 
\end{eqnarray}
}}

\noindent where $L_i = (l_1,\ldots,l_i)\in (\NN^{\ast})^i$, with  $l_1+\ldots +l_i = j$, is an $i$-dimensional multi-index and 
{\footnotesize {
\begin{eqnarray}\nonumber  && \hspace{-0.95cm} \displaystyle deg_{L_i}\left[D_1,\ldots, D_r\right] =\\\nonumber  &=& \hspace{-0.35cm}\left(\displaystyle\sum_{1\leq e_{1}^1<\ldots<e_{l_1}^1 \leq r}deg(D_{e_{1}^1}) \cdot\ldots\cdot deg(D_{e_{l_1}^1})\right)\hspace{-0.12cm}\cdot \ldots \cdot \hspace{-0.12cm} \left(\displaystyle\sum_{1\leq e_{1}^i<\ldots<e_{l_i}^i \leq r}deg(D_{e_{1}^i}) \cdot \ldots \cdot deg(D_{e_{l_i}^i})\hspace{-0.12cm}\right).
\end{eqnarray}
}}
\end{cor}

\noindent {\bf Proof:}  For each $\lambda = 1,\ldots, r$, we have
$$
[D_{\lambda}] = \OO_{\PP^m}(1)^{\otimes deg(D_{\lambda})}
$$ 
\noindent where $\OO_{\PP^m}(1)$ denotes the hyperplane bundle. Consequently,
$$
c_1([D_{\lambda}]) = deg\,(D_{\lambda})c_1(\OO_{\PP^m}(1))
$$
\noindent and
$$
\prod_{\lambda = 1}^rc_1([D_{\lambda}]) = \prod_{\lambda = 1}^rdeg(D_{\lambda})c_1(\OO_{\PP^m}(1))^{r}.
$$
\noindent Thus it follows from Corollary \ref{coroo} that 
\begin{eqnarray} \nonumber
 \hspace{-1.9cm}&& \hspace{-0.9cm}\displaystyle GSV(\fol,\hspace{-0.05cm}V) = \hspace{-0.1cm}\int_{\PP^m} \hspace{-0.09cm} \prod_{\lambda = 1}^r\hspace{-0.09cm}deg(D_{\lambda})c_1(\OO_{\PP^m}(1))^{r}\hspace{-0.08cm}\sum_{t=0}^{m-r}\hspace{-0.07cm} ( \hspace{-0.01cm} c_t(\PP^n) \hspace{-0.05cm} + \\\nonumber&& \hspace{-0.9cm} + \hspace{-0.05cm}\sum_{j=1}^t\sum_{i=1}^j\sum_{\mid L\mid=j}(-1)^ic_{t-j}(\PP^n)c_{L_i}(\bigoplus^r_{\lambda=1}[D_{\lambda}]))c_1(T\fol^{\ast})^{m-r-t} = \\ \nonumber  && \hspace{-0.9cm} =  \int_{\PP^m}\hspace{-0.09cm} \prod_{\lambda = 1}^r\hspace{-0.09cm}deg(D_{\lambda})\hspace{-0.15cm}\sum_{t=0}^{m-r}( {m+1 \choose t}c_1(\OO_{\PP^m}(1))^{t} \hspace{-0.1cm} + \\\nonumber && \hspace{-1cm}  + \hspace{-0.1cm}\sum_{j=1}^t\sum_{i=1}^j\sum_{\mid L_i\mid=j}\hspace{-0.15cm}(-1)^i\hspace{-0.07cm}{m+1 \choose t-j}c_1(\OO_{\PP^m}(1))^{t-j}c_{L_i}(\bigoplus^r_{\lambda=1}[D_{\lambda}]))(d-1)^{m-r-t}\hspace{-0.05cm}c_1(\OO_{\PP^m}(1))^{m-t}
\end{eqnarray}

\noindent where in last equality we use the relations $$c_1(T\fol^{\ast}) = (d-1)c_1(\OO_{\PP^m}(1))$$ and $$c_{t-j}(\PP^m) = {m+1 \choose t-j}c_1(\OO_{\PP^m}(1))^{t-j},\,\,\, \mbox{for all $j=0,\ldots,t$.}$$

\noindent On the other hand, for all $L_i = (l_1,\ldots,l_i)\in (\NN^{\ast})^i$, with  $l_1+\ldots +l_i = j$, we have that {\small{ 
\begin{eqnarray}\nonumber && \hspace{-0.7cm}c_{L_i}\left(\bigoplus^r_{\lambda=1}[D_{\lambda}]\right) = c_{l_1}\left(\bigoplus^r_{\lambda=1}[D_{\lambda}]\right)\cdots c_{l_i}\left(\bigoplus^r_{\lambda=1}[D_{\lambda}]\right) = \\\nonumber && \hspace{-0.7cm} = \left(\displaystyle\sum_{1\leq e_{1}^1<\ldots<e_{l_1}^1 \leq r}\hspace{-0.7cm}c_1([D_{e_{1}^1}]) \ldots c_1([D_{e_{l_1}^1}])\right) \ldots  \left(\displaystyle\sum_{1\leq e_{1}^i<\ldots<e_{l_i}^i \leq r}\hspace{-0.7cm}c_1([D_{e_{1}^i}])  \ldots  c_1([D_{e_{l_i}^i}])\right) = \\\nonumber && \hspace{-0.7cm} = \left(\displaystyle\sum_{1\leq e_{1}^1<\ldots<e_{l_1}^1 \leq r} \hspace{-0.7cm} deg(D_{e_{1}^1}) \ldots deg(D_{e_{l_1}^1})c_1(\OO_{\PP^m}(1))^{l_1}\right)\ldots\left(\displaystyle\sum_{1\leq e_{1}^i<\ldots<e_{l_i}^i \leq r}\hspace{-0.7cm} deg(D_{e_{1}^i})  \ldots  deg(D_{e_{l_i}^i})c_1(\OO_{\PP^m}(1))^{l_i})\right) = \\\nonumber && \hspace{-0.7cm} = \hspace{-0.1cm} \left(\hspace{-0.05cm}\displaystyle\sum_{1\leq e_{1}^1<\ldots<e_{l_1}^1 \leq r}\hspace{-0.7cm}deg(D_{e_{1}^1})  \ldots deg(D_{e_{l_1}^1})\hspace{-0.12cm}\right) \hspace{-0.10cm}\ldots \hspace{-0.10cm}  \left(\hspace{-0.05cm}\displaystyle\sum_{1\leq e_{1}^i<\ldots<e_{l_i}^i \leq r}\hspace{-0.6cm} deg(D_{e_{1}^i}) \ldots  deg(D_{e_{l_i}^i})\right)c_1(\OO_{\PP^m}(1))^{j} 
\end{eqnarray}}}

\noindent and, we get $c_{L_i}\left(\bigoplus^r_{\lambda=1}[D_{\lambda}]\right) = deg_{L_i}\,[D_1,\ldots, D_r]c_1(\OO_{\PP^m}(1))^{j}$. Therefore, we obtain
{\small{
\begin{eqnarray} \nonumber
& &\hspace{-0.7cm}  GSV(\fol,V) = \\\nonumber & & \hspace{-0.7cm} = \hspace{-0.08cm} \prod_{\lambda = 1}^r\hspace{-0.07cm}deg(D_{\lambda})\hspace{-0.15cm}\sum_{t=0}^{m-r}\hspace{-0.15cm}\left( \hspace{-0.15cm} {m+1 \choose t} \hspace{-0.08cm} + \hspace{-0.08cm} \sum_{j=1}^t\hspace{-0.08cm} \sum_{i=1}^j \hspace{-0.11cm}\sum_{\mid L\mid=j}\hspace{-0.1cm}(-1)^i\hspace{-0.04cm}{m+1 \choose t-j}c_{L_i}\hspace{-0.1cm}\left(\bigoplus^r_{\lambda=1}[D_{\lambda}]\hspace{-0.1cm} \right)\hspace{-0.15cm} \right)\hspace{-0.1cm}(d-1)^{m-r-t}\hspace{-0.2cm}\int_{\PP^m}\hspace{-0.3cm}c_1(\OO_{\PP^m}(1))^{m}  
\end{eqnarray}

}}

\noindent and the proof follows from the fact that $\int_{\PP^m}c_1(\OO_{\PP^m}(1))^m =1$.  

\noindent $\square$

\begin{exe}
Consider the vector field $v$ in $\mathbb{C}^4$ defined by
\[
v = z_0 \frac{\partial}{\partial z_0} + 7z_1 \frac{\partial}{\partial z_1} + 3z_2 \frac{\partial}{\partial z_2} + 4z_3 \frac{\partial}{\partial z_3}.
\]
This vector field induces a one-dimensional holomorphic foliation $\mathcal{F}$ on $\mathbb{P}^3$ of degree $d = 1$. The singular set of $\mathcal{F}$ is given by
\[
\operatorname{Sing}(\mathcal{F}) = \{[1:0:0:0],\quad [0:1:0:0],\quad [0:0:1:0],\quad [0:0:0:1]\}.
\]

\noindent Let $C \subset \mathbb{P}^3$ be the $\mathcal{F}$-invariant complete intersection curve defined by
\[
C = \{ [z_0:z_1:z_2:z_3] \in \mathbb{P}^3 \mid z_0^2 z_1 - z_2^3 = 0 \quad \text{and} \quad z_3^2 - z_0 z_1 = 0 \},
\]
of multi-degree $(k_1,k_2) = (3,2)$. The singular set of $C$ is
\[
\operatorname{Sing}(C) = \{[1:0:0:0],\quad [0:1:0:0]\}.
\]

\noindent We compute the total GSV index of the foliation $\mathcal{F}$ with respect to $C$:
\[
GSV(\mathcal{F}, C) = GSV_{[1:0:0:0]}(\mathcal{F},C) + GSV_{[0:1:0:0]}(\mathcal{F},C).
\]

\noindent In the affine chart $z_0\neq 0$ with coordinates $(x_1,x_2,x_3) = (z_1/z_0, z_2/z_0, z_3/z_0)$, the foliation is defined by
\[
v_{0} = 6x_1 \frac{\partial}{\partial x_1} + 2x_2 \frac{\partial}{\partial x_2} + 3x_3 \frac{\partial}{\partial x_3},
\]
and the curve $C$ is defined by
\[
C_0 = \{(x_1,x_2,x_3) \in \mathbb{C}^3 : x_1 - x_2^3 = 0 \quad \text{and} \quad x_3^2 - x_1 = 0\}.
\]
Following \cite[Theorem 2.4 and Corollary 4.2]{aaaa}, we obtain:
{\small{\begin{eqnarray}\nonumber 
& & \hspace{-0.65cm} GSV_{[1:0:0:0]}(\fol,C) =\\\nonumber 
&& \hspace{-0.65cm} = - \dim_{\CC}\frac{\OO_{3,0}}{\left\langle -3x_2^2, 2x_3, -6x_3x_2^2, x_1 - x_2^3, x_3^2 - x_1 \right\rangle} + \dim_{\CC}\frac{\OO_{3,0}}{\left\langle 6x_1, 2x_2, 3x_3, x_1 - x_2^3, x_3^2 - x_1 \right\rangle}
\\\nonumber 
&& \hspace{-0.65cm} = -2 + 1 = -1
\end{eqnarray}}}

\noindent In the affine chart $z_1\neq 0$ with coordinates $(y_1,y_2,y_3) = (z_0/z_1, z_2/z_1, z_3/z_1)$, $\mathcal{F}$ is defined by
\[
v_1 = -6y_1 \frac{\partial}{\partial y_1} - 4y_2 \frac{\partial}{\partial y_2} - 3y_3 \frac{\partial}{\partial y_3},
\]
and $C$ is given by
\[
C_1 = \{(y_1,y_2,y_3) \in \mathbb{C}^3 : y_1^2 - y_2^3 = 0 \quad \text{and} \quad y_3^2 - y_1 = 0\}.
\]
In this case, for the point $[0:1:0:0]$, we have:
{\small{\begin{eqnarray}\nonumber 
& & \hspace{-0.65cm} GSV_{[0:1:0:0]}(\fol,C) =\\\nonumber 
&& \hspace{-0.65cm} = - \dim_{\CC}\frac{\OO_{3,0}}{\left\langle -3y_2^2, 4y_1y_3, -6y_2^2y_3, y_1^2 - y_2^3, y_3^2 - y_1 \right\rangle} + \dim_{\CC}\frac{\OO_{3,0}}{\left\langle -6y_1, - 4y_2, -3y_3, y_1^2 - y_2^3, y_3^2 - y_1 \right\rangle}
\\\nonumber 
&& \hspace{-0.65cm} = -6 + 1 = -5.
\end{eqnarray}}}

\noindent Consequently, we obtain the total $GSV(\mathcal{F},C) = -6$. We used the software \textit{Singular} \cite{singular} to compute the local indices $GSV_{[1:0:0:0]}(\mathcal{F},C)$ and $GSV_{[0:1:0:0]}(\mathcal{F},C)$.

\noindent We observe that the formula provided in Corollary \ref{GSformula} offers a more direct method for calculating $GSV(\mathcal{F},C)$:
\[
GSV(\mathcal{F},C) = (k_1 k_{2})\left[(m+1) - (k_1 + k_{2}) +(d-1) \right] = -6.
\]
Furthermore, in this case, we have $GSV < 0$ and $k_1 + k_2 > d+m$, in accordance with Theorem 3.
\end{exe}

\section{Bounding the GSV index}\label{secc03}

In this section, we present the proof of Theorem 2, which provides bounds for the GSV index of the foliation. For this, let us consider a one-dimensional holomorphic foliation $\fol$, with isolated singularities, leaving invariant a complete intersection
$V$, of codimension $r$, with isolated singularities. 

Given $p\in Sing(\fol) \cap Sing(V)$, we consider $f = (f_1,\ldots, f_r) = 0$ a local reduced equation of $V$ on a neighborhood $U$ of $p$ in $X$, with a local coordinate system $(z_1,\ldots,z_m)$,  and (by taking a smaller $U$, if necessary) $v = \sum_{i=1}^ma_i\frac{\partial}{\partial z_i}$ a holomorphic vector field defining the foliation $\fol$ on $U$. Since $V$ is invariant by $\fol$, there exists a $r\times r$ matrix $H = (h_{ij})_{r\times r}$ with holomorphic entries such that

$$
df(v) = \left\langle H, f \right\rangle.
$$

\noindent  According to \cite[Theorem 2.4 and Corollary 4.2] {aaaa}, GSV index of $\fol$ in $p$ along $V$ can be calculated by the following algebraic formula

{\small{
\begin{eqnarray} \label{revrev01}
 \,\,\,\,\,\,\,\,\,\,\,\,GSV_p(\fol,V) = \varepsilon_r\, \dim \frac{\OO_{m,p}}{\left\langle v \right\rangle}
+(-1)^{m-r}  \tau_p(V) +(-1)^{m-r-1} \dim\, \mbox{coker} (\gamma_{m-r-1}),
\end{eqnarray}}}

\medskip

\noindent where $\displaystyle \varepsilon_r$ denotes the integer

$$
\displaystyle \varepsilon_r = \left\{
\begin{array}{lll}
\displaystyle\sum_{j=0}^{m-r-2} (-1)^j {r-1+j \choose j},\,\,\mbox{if $r \leq m - 1$}\\
0,\,\,\mbox{if $r = m - 1$,}
\end{array} \right.
$$

\medskip

\noindent $\OO_{m,p}$ denotes the ring of germs of holomorphic functions at $(X,p)$,  $\left\langle v \right\rangle \subset \OO_{m,p}$ is the ideal generated by the components $a_1,\ldots, a_m$ of vector field $v$, $\tau_p(V)$ is the invariant of Greuel \cite{ggel},  which can be calculated by

$$
\tau_p(V) = \displaystyle dim_{\CC}\frac{\OO_{m,p}}{\left\langle f,\displaystyle\frac{\partial(f_1,\ldots,f_r)}{\partial(x_{j_1},\ldots,x_{j_r})} (1\leq j_1<\ldots<j_r \leq m)\right\rangle}
$$

\noindent and $\gamma_{m-r-1}$ is the $(m-r-1)$-st map of {\it small Gobelin complex} (see \cite{aaaa})

{\footnotesize{
$$
$$
\[
0
\xleftarrow{}
\displaystyle\left(\frac{\OO_{m,p}}{\left\langle v \right\rangle}\right)^{{r-1 \choose 0}}
\xleftarrow[\gamma_0]{\mathbb{M}(0)= \left  (\begin{smallmatrix}{{f}}_{1}&
      {{f}}_{{2}}&
	\cdots &
      {{f}}_{{r}}\\
      \end{smallmatrix}\right)}
\left(\frac{\OO_{m,p}}{\left\langle v \right\rangle}\right)^{{r + 0 \choose 1}}
\xleftarrow[\gamma_1]{  \mathbb{M}(1)= \left(\begin{smallmatrix}{{f}}_{{2}}& 0 & f_3 &
      {\cdots}&
      0& 0 &
      {{f}}_{{r}}&
      {{c}}_{1,1}&
      {\cdots}&
      {{c}}_{{r},1}\\
      {-{{f}}_{1}}& f_3  & 0 &
      {\cdots}&
      0 & {{f}}_{{r}}&
      0&
      {{c}}_{1,{2}}&
      {\cdots}&
      {{c}}_{{r},{2}}\\
       0 &  -f_2& -f_1  &
      {\cdots}&
      f_r & 0&
      0&
      {{c}}_{1,{3}}&
      {\cdots}&
      {{c}}_{{r},{3}}\\
    \cdot & \cdot & \cdot & & \cdot & \cdot & \cdot & \cdot & & \cdot\\
      \cdot & \cdot & \cdot & & \cdot & \cdot & \cdot & \cdot & & \cdot\\
      0& 0 & 0 &
      {\cdots}&
      -f_3 & {-{{f}}_{{2}}}&
      {-{{f}}_{1}}&
      {{c}}_{1,{r}}&
      {\cdots}&
      {{c}}_{{r},{r}}\\
      \end{smallmatrix}\right) }\left(\frac{\OO_{m,p}}{\left\langle v \right\rangle}\right)^{{r+1 \choose 2}}   \ldots\, .
\]
$$
$$
}}

\noindent We recall that each map $\gamma_i$ is defined by a matrix product

\begin{eqnarray}\nonumber
\gamma_i: \displaystyle\left(\frac{\OO_{m,p}}{\left\langle v \right\rangle}\right)^{{r + i \choose i+1}} \displaystyle &\longrightarrow& \displaystyle\left(\frac{\OO_{m,p}}{\left\langle v \right\rangle}\right)^{{r + i -1 \choose i}}\\\nonumber \\\nonumber
g  &\longmapsto& \gamma_i (g) = [{\mathbb{M}(i)}][g]^T
\end{eqnarray}
\noindent  with $[{\mathbb{M}(i)}]$ being an appropriate matrix of order ${r+i-1 \choose i} \times {r+i \choose i+1}$ whose non-zero entries are formed by the components of $f=(f_1,\ldots,f_r)$ and by the entries of matrix $H = (h_{ij})_{r\times r}$ (which is defined by relation $df(v) = \left\langle H, f \right\rangle$). 

Furthermore, if for each row $R_j$ $\left(1\leq  j \leq {r+i-1 \choose i}\right)$ of the matrix ${\mathbb{M}(i)}$ we denote by $\left\langle R_j \right\rangle$ the ideal generated by its elements, then we can write

$$
\mbox{coker} (\gamma_{m-r-1}) = \prod_{j=1}^{{m-2 \choose m-r-1}} \frac{\OO_{m,p}}{\left\langle v, R_j \right\rangle}
$$ 
\noindent and
$$
\dim\, \mbox{coker} (\gamma_{m-r-1}) = \sum_{j=1}^{{m-2 \choose m-r-1}} \dim_{\CC}\frac{\OO_{m,p}}{\left\langle v, R_j \right\rangle}.
$$ 
\noindent Thus,  considering the integer $\rho: = \#\{\left\langle R_j \right\rangle: \left\langle R_j \right\rangle \subset \left\langle v \right\rangle  \}$ we get that
\begin{eqnarray} \label{revrev02} 
\displaystyle \dim\, \mbox{coker} (\gamma_{m-r-1})  =  \rho\,\dim_{\CC}\frac{\OO_{m,p}}{\left\langle v \right\rangle}  + \sum_{\left\langle R_j \right\rangle \not\subset \left\langle v \right\rangle } \dim_{\CC}\frac{\OO_{m,p}}{\left\langle v, R_j \right\rangle}. 
\end{eqnarray}



\begin{remark} In the case $V = C$ is a complete intersection curve ($codim(V) = m-1$), formula (\ref{revrev01}) reduces to  
$$
GSV_p(\fol, C) \,\,\,=\,\,\, -\, \tau_p(C) + \dim\, \mbox{coker} (\gamma_{0}),
$$
\noindent where the map $\gamma_0$ is given by 

$$
\displaystyle \frac{\OO_{m,p}}{\left\langle v \right\rangle}
\xleftarrow[\gamma_0]{\left  (\begin{smallmatrix}{{f}}_{1}&
      {{f}}_{{2}}&
	\cdots &
      {{f}}_{{m-1}}\\
      \end{smallmatrix}\right)}
\left(\frac{\OO_{m,p}}{\left\langle v \right\rangle}\right)^{m-1}
$$

\noindent and therefore, we obtain
\begin{eqnarray}\label{re002}
GSV_p(\fol, C) \,\,\,=\,\,\, -\, \tau_p(C) + \dim_{\CC}\frac{\OO_{m,p}}{\left\langle v, f \right\rangle}
\end{eqnarray}
\end{remark}

\medskip

\noindent {\bf Proof of Theorem 2:} Suppose that the singularities of $\fol$ are non-degenerate. It follows from (\ref{revrev01}) and (\ref{revrev02}) that
\begin{eqnarray} \nonumber
GSV_p(\fol,V) &=& \varepsilon_r
+(-1)^{m-r}  \tau_p(V) +(-1)^{m-r-1} \left(\rho  + \sum_{\left\langle R_j \right\rangle \not\subset \left\langle v \right\rangle } \dim_{\CC}\frac{\OO_{m,p}}{\left\langle v, R_j \right\rangle} \right)\\\nonumber &=& \varepsilon_r
+(-1)^{m-r}  \tau_p(V) +(-1)^{m-r-1} \rho
\end{eqnarray}
\noindent where in the last step we have used the fact that
$$
\left\langle R_j \right\rangle \not\subset \left\langle v \right\rangle \Longrightarrow \dim_{\CC}\frac{\OO_{m,p}}{\left\langle v, R_j \right\rangle} \,\,\,\leq \,\,\,\dim_{\CC}\frac{\OO_{m,p}}{\left\langle v \right\rangle} - 1 \,\,\,=\,\,\, 0.
$$
\noindent If $m-r$ is even, then
\begin{eqnarray}\label{revrev04}
GSV_p(\fol,V) = \varepsilon_r
+  \tau_p(V) - \rho.
\end{eqnarray}

\noindent Since $\rho\in \left\{0,1,\ldots, {n-2 \choose n-r-1}\right\}$, we obtain

{{\small $$
\tau_p(V) + \sum_{j=0}^{m-r-1} (-1)^j {r-1+j \choose j} \,\,\, \leq \,\,\, GSV_p(\fol, V) \,\,\, \leq \,\,\, \tau_p(V) + \sum_{j=0}^{m-r-2} (-1)^j {r-1+j \choose j}.
$$}} 

\noindent If $m-r$ is odd, then
\begin{eqnarray}\label{revrev05}
GSV_p(\fol,V) = \varepsilon_r
-  \tau_p(V) + \rho
\end{eqnarray}

\noindent and since $\rho\in \left\{0,1,\ldots, {m-2 \choose m-r-1}\right\}$, we obtain

{{\small $$
\sum_{j=0}^{m-r-2} (-1)^j {r-1+j \choose j} - \tau_p(V) \,\,\, \leq \,\,\, GSV_p(\fol, V) \,\,\, \leq \,\,\, \sum_{j=0}^{m-r-1} (-1)^j {r-1+j \choose j} - \tau_p(V).
$$}} 

\noindent $\square$
{\small{

\begin{table}[htb]
\hspace*{1.375cm } \begin{tabular}{|c|c|}
      \hline
      \,\,\,\,\,\,\,\,\,\,\,\,\,\,\,\hspace{1.0cm } {\bf lower bound} \hspace{0.983cm }\,\,\,\,\,\,\,\,\,\,\,\,\,\, & \,\,\,\,\,\,\,\,\,\,\,\,\,\,\,\hspace{0.55cm }{\bf upper bound}\hspace{0.55cm }\,\,\,\,\,\,\,\,\,\,\,\,\, \\
			       \end{tabular}
\vspace{-0.67cm}
\begin{center}
\begin{tabular}{ | c | c | c | c | c |  c |} 
\hline
$dim(V)$ & $m$ even & $m$ odd & $m$ even & $m$ odd \\ 
\hline
$m-1$ & $(m-1) + \tau_p(V)$ & $(m-2) - \tau_p(V)$ & $m + \tau_p(V)$ & $(m-1) - \tau_p(V)$\\ \hline
$m-2$ & $-\frac{(m-2)}{2} + \tau_p(V)$ & $-\frac{(m-1)}{2} +1 - \tau_p(V)$ & $\frac{(m-2)}{2} + \tau_p(V)$ & $\frac{(m-1)}{2} - \tau_p(V)$\\ \hline
$2$ & $-(m-3)+\tau_p(V)$ & $-(m-3)+\tau_p(V)$ & $1+\tau_p(V)$  & $1+\tau_p(V)$\\ \hline
$1$ &  $- \tau_p(V)$ & $- \tau_p(V)$ & $1- \tau_p(V)$ & $1 - \tau_p(V)$\\ \hline
\end{tabular}
\end{center}
\caption{A summary of the bounds of the $GSV_p(\fol,V)$,  according to Theorem 2.}
\label{xis}
\end{table}

}}

\begin{remark} \label{rr1} 
\noindent{\bf (a)} According to item (ii) of Theorem 2, if $dim(V) =1$ then 
$$
-\tau_p(V) \leq GSV_p(\fol,V) \leq 1 -\tau_p(V).
$$ 
\noindent In particular, $GSV_p(\fol,V) \leq 0$.  Table \ref{xis} provides a summary of the bounds of the GSV index for certain cases, according to Theorem 2.\\\\
\noindent{\bf (b)} In the Theorem 2, if $m-r$ is even (resp., if $m-r$ is odd) the GSV index achieves upper bound (resp., lower bound) if $\rho = 0$, i.e., at least one element in each row of the matrix of map $\gamma_{m-r-1}$ does not belong to the ideal $\left\langle v \right\rangle$. In addition, if $m-r$ is even (resp., if $m-r$ is odd) the GSV index achieves lower bound (resp., upper bound) if $\rho = {m-2 \choose m-r-1}$, i.e., all the elements of the matrix of map $\gamma_{m-k-1}$ belong to the ideal  $\left\langle v \right\rangle$.
\end{remark}

Relations (\ref{revrev04}) and (\ref{revrev05}) give us the following characterization of the positivity of the GSV index:

\begin{cor}\label{cor0002}
Let $\fol$ be a one-dimensional holomorphic foliation, with nondegenerate isolated singularities, over an $m$-dimensional complex manifold $X$, and $V\subset X$ a local complete intersection  of codimension $r$, with isolated singularities, invariant by $\fol$.  For all $p\in Sing(\fol) \cap Sing(V)$ one has

\medskip

\noindent (i) If $m-r$ is even then
$$
GSV_p(\fol, V)> 0 \Longleftrightarrow
\tau_p(V) + \varepsilon_r> \rho
$$

\bigskip

\noindent (ii) If $m-r$ is odd then
$$
GSV_p(\fol, V)> 0 \Longleftrightarrow
\tau_p(V) - \varepsilon_r < \rho
$$
\end{cor}

\section{Applications}\label{secc04}

In this section, we turn to projective foliations with the aim of obtaining applications in the context of the Poincaré Problem. Our goal is to present the proofs of Theorem 3, Corollary 1.3, and Theorem 4. Moreover, we shall present a few remarks  on Schwartz index in the context of holomorphic foliations.

Let $\fol$ be a one-dimensional foliation on $\PP^n$, with isolated singularities and degree $\deg(\fol) = d$. If $C \subset \PP^m$  is an $\fol$ - invariant complete intersection curve of multi-degree $(k_1,\ldots, k_{m-1})$, with isolated
singularities, it follows from Corollary \ref{GSformula} that 
\begin{eqnarray} \nonumber
GSV(\fol,V) = (k_1\cdots k_{m-1})\left[(m+1) - \left(k_1 +\cdots + k_{m-1}\right) +(d-1) \right].
\end{eqnarray}

\noindent Thus, we get $GSV(\fol,V)\geq 0$ if and only if 
$\displaystyle k_1 +\cdots + k_{m-1} \leq d + m$, according to Theorem 3.

We observe that if $C$ is smooth, then at each singular point of $\fol$ the GSV index reduces to the (classical) Poincaré--Hopf index of the local vector field defining $\fol$, restricted to $C$, which is positive. Hence, we obtain $GSV(\fol, V) \ge 0$ and, by Theorem 3, the conditions $\displaystyle k_1 +\cdots + k_{m-1} \leq d + m$ occur, according to Corollary 1.3.

Now, to prove Theorem 4, we will assume that $\fol$ has nondegenerate isolated singularities. By using formula (\ref{re002}), we obtain
\begin{eqnarray}\nonumber
GSV_p(\fol, C) \,\,\,=\,\,\, -\, \tau_p(C) + \dim_{\CC}\frac{\OO_{m,p}}{\left\langle v, f \right\rangle} \,\,\,\geq\,\,\, -\, \mu_p(C) + 1
\end{eqnarray}
\noindent where in the last step we have used the fact that $\mu_p(C) \geq \tau_p(C)$ (see \cite{ggel}) and $\dim_{\CC}\frac{\OO_{m,p}}{\left\langle v, f \right\rangle} \geq 1$. Thus, we get
$$
 (k_1\cdots k_{m-1})\left[(m+1) - \left(k_1 +\cdots + k_{m-1}\right) +(d-1) \right] 
\geq \,\,\, -\, \sum_{p\in Sing(C)}(\mu_p(C) - 1)
$$
\noindent and, consequently,
$$
(k_1 \ldots k_{m-1}) \left(k_1 +\ldots + k_{m-1} - m \right)  \hspace{0.0cm} - \sum_{p\in Sing(C)}\hspace{-0.35cm}(\mu_p(C) -1) \leq d(k_1 \ldots k_{m-1}),
$$
\noindent according to Theorem 4.

\begin{remark}
We observe that the inequality $\displaystyle k_1 +\cdots + k_{m-1} \leq d + m$ is similar to those in \cite{cja}, \cite{cm3} and \cite{CE}.
\end{remark}

\subsection{Some remarks on Schwartz index}\label{sch01}

M. -H. Schwartz carried out pioneering work (\cite{SSS1}, \cite{SSS2}) on indices of vector fields defined on singular varieties, introducing the notion of the radial index for a special class of vector fields and establishing a version of the Poincaré-Hopf theorem for singular varieties. This notion was subsequently generalized by various authors in different settings \cite{AGG}, \cite{ABC00}, \cite{ABC}, \cite{Ebe}, \cite{Ebe2}, \cite{King}, \cite{SeaSuw2}, and is also referred to in the literature as the Schwartz index. Analogously to the GSV index, the Schwartz index is well suited to the study of singular holomorphic foliations and is used to study its residues relative to singular invariant subvarieties (see for example \cite{BarSaeSuw}, \cite{DSM2} and \cite{Suw2}).

Let $X$ be a complex manifold and $\fol$ a one-dimensional holomorphic foliation on $X$, with isolated singularities, leaving invariant a variety $V$, with isolated singularities. For each $p\in V$ we denotes by $Sch_p(\fol,V)$ the {\it Schwartz index of $\fol$ at $p$ along $V$}. For definition and details on the Schwartz index in the context of foliations we refer to \cite[Ch. IV, Sec. 7] {Suw2} 

We know that $Sch_p(\fol, V)$ is positive when $\fol$ is locally defined by radial vector fields (see \cite[Ch.~IV, Definition~1.3]{Suw2}) or when $V$ is a hypersurface and $\dim(X)$ is even (see \cite[Theorem~2]{DSM2}). However, it remains unknown whether this positivity holds in full generality. Theorem~5 asserts that the index is positive when the invariant variety is a complete intersection curve:

\begin{teo5}\label{hhh}
Let $X$ be an $m$-dimensional complex manifold, $C\subset X$ a local complete intersection curve and $\fol$ a one-dimensional holomorphic foliation on $X$, with isolated singularities, which leaves $C$ invariant. One has 
$$Sch_p(\fol, C)>0, \,\,\mbox{for all}\,\, p\in Sing(\fol).$$

In particular, when $X = \PP^m$, $Sch_p(\fol, C)\geq 2$, if $(C,p)$ is not quasi-homogeneous germ.
\end{teo5}

\noindent {\bf Proof:} Let $f = (f_1,\ldots, f_{m-1}) = 0$ be a local reduced equation of $C$ on a neighborhood $U$ of $p$ in $X$, with a local coordinate system $(z_1,\ldots,z_m)$,  and $v = \sum_{i=1}^ma_i\frac{\partial}{\partial z_i}$ a holomorphic vector field defining the foliation $\fol$ on $U$. According to relation (\ref{re002}) we have
$$
GSV_p(\fol, C) \,\,\,=\,\,\, -\, \tau_p(C) + \dim_{\CC}\frac{\OO_{m,p}}{\left\langle v, f \right\rangle}.
$$
\noindent By \cite[Ch.IV, Proposition 1.9]{Suw2},
$$
GSV_p(\fol,C) = Sch_p(\fol,C)  - \mu_p(C),
$$
\noindent where $\mu_p(C)$ denotes the Milnor number of $C$ at $p$. Since $\mu_p(C) \geq \tau_p(C)$ (see \cite{ggel}) and $\dim_{\CC}\frac{\OO_{m,p}}{\left\langle v, f \right\rangle} \geq 1$, we obtain
\begin{eqnarray}\nonumber
Sch_p(\fol,C) &=& \mu_p(C) - \tau_p(C) + \dim_{\CC}\frac{\OO_{m,p}}{\left\langle v, f \right\rangle} \,\,\,> \,\,\,0 .
\end{eqnarray}

On the other hand, if $X = \PP^m$, we know that $\mu_p(C) = \tau_p(C)$ if and only if  $(C,p)$ is quasi-homogeneous (see \cite{henrik}). Hence, if $(C,p)$ is not quasi-homogeneous, then
\begin{eqnarray}\nonumber
Sch_p(\fol,C) = \mu_p(C) - \tau_p(C) + \dim_{\CC}\frac{\OO_{m,p}}{\left\langle v, f \right\rangle} \,\,\,\geq \,\,\, 2.
\end{eqnarray}

\noindent $\square$

As a consequence of Theorem 5, we will show that the Euler characteristic provides a specific obstruction to the existence of holomorphic vector fields on a complete intersection curve.

\begin{cor}
Let $C$ be a complete intersection curve with isolated singularities in a complex manifold. If $C$ admits a holomorphic vector field with $l$ singularities, then
$$
\chi(C) \geq l.
$$
\end{cor}

\noindent {\bf Proof:} Suppose that $C$ admits a holomorphic vector field $v$ with $Sing(v) = \{p_1,\ldots, p_l\}$. By using \cite[Theorem 2.1.1]{BarSaeSuw}, we have
$$
\chi (C) = \sum_{j=1}^{l} Sch_{p_j}(v, C).
$$
\noindent Since $Sch_{p_j}(\fol, C)>0, \,\,\mbox{for all}\,\, p_j\in Sing(\fol),$ we obtain
$$
\chi (C) \geq l.
$$
\noindent $\square$

 \begin{remark} In \cite{DSM2} we prove that the obstruction determined by the Euler characteristic for the existence of vector fields is even more comprehensive in the case of hypersurfaces with isolated singularities.
\end{remark}

{\it Acknowledgements.} I am grateful to D. Lehmann for interesting conversations.


\begin{thebibliography}{99}

\bibitem{AGG} Aguilar, M.,  Seade J.  and Verjovsky A.: Indices of vector fields and topological invariants of real analytic singularities, Crelle's, J. Reine u. Ange. Math. 504
(1998), 159-176.


\bibitem{BB1} Baum, P.  and  Bott R.: Singularities of Holomorphic Foliations, J. Differential Geom, 7, 279-342, 1972.


\bibitem{ABC00} Brasselet, J.-P., Lehmann, D., Seade, J. and Suwa, T.: Milnor numbers and classes of local complete intersections, Proc. Japan Acad. Ser. A Math. Sci 75 (1999),
179-183.


\bibitem{ABC} Brasselet, J.-P., Lehmann, D., Seade, J. and Suwa, T.: Milnor classes of local complete intersections, Transactions A. M. S. 354 (2001), 1351-1371.


\bibitem{a03} Brasselet, J.-P., Seade, J. and Suwa, T.: An explicit cycle representing the Fulton-Johnson class, Singularit\'es Franco-Japonaises, S\'emin. Congr., 10, Soc. Math. France, Paris, p. 21-38, 2005.

\bibitem{BarSaeSuw} Brasselet, J.-P.,  Seade, J.  and   Suwa T.: Vector Fields on Singular Varieties, Lecture Notes in Mathematics, Spring, 2009.

\bibitem{Bru1} Brunella, M.: Feuilletages holomorphes sur les surfaces complexes compactes. Ann. Sci. E.N.S 30 (1997), 569-594.

\bibitem{Bru} Brunella, M.: Some remarks on indices of holomorphic vector fields. Publ. Math. Vol. 41, No. 2, 527-544, 1997.

\bibitem{a09} Brunella, M.: Birational Geometry of Foliations, Publica\c c\~oes Matem\'aticas, IMPA, Rio de Janeiro, 2010.


\bibitem{aaaa}  Bothmer, H. C. G., Ebeling, W. and G\'omez-Mont, X.: An algebraic formula for the index of a vector field on an isolated complete intersection singularity, J. Algebraic Geom. 7, 731-752, 1998.

\bibitem{Can} Carnicer, M. M.: The Poincar\'e problem in the nondicritical case. Ann. Math. Vol. 140, No. 2, 289-294, 1994.

\bibitem{CL} Cavalier, V. and Lehmann, D.: On the Poincar'e inequality for one-dimensional foliations. Compositio Math. 142 (2006) 529-540.

\bibitem{Alc} Cerveau, D.  and Lins Neto, A.: Holomorphic foliations in $\mathbb{P}^2$ having an invariant algebraic curve. 
Ann. Inst. Fourier. Vol. 41, No. 4, 883-903, 1991.

\bibitem{cja} Corr\^ea, M. and  Jardim, M.: Bounds for sectional genera of varieties invariant under Pfaff fields,  Illinois Journal of Mathematics, 
vol. \textbf{56} (2) (2012),   343-352.


\bibitem{cm3} Corr\^ea, M. and  Machado,  D.: GSV-Index for Holomorphic Pfaff Systems. Doc. Math. 25 (2020), pp. 1011-1027.


\bibitem{MD3} Corr\^ea, M. and  Machado,  D.: A global residue formula for logarithmic indices of  foliations. To appear in Communications in Analysis and Geometry, Vol. 32, No. 4, 2024.


\bibitem{CE}  Cruz, J.  D. A. and  Esteves, E.: Regularity of subschemes invariant under Pfaff fields on projective spaces. Commentarii Mathematici Helvetici \textbf{ 86} (2011),  947-965. 


\bibitem{singular}  Decker, W., Greuel, G.-M. , Pfister, G.  and Schönemann H., Singular 4-4-0 - A computer algebra system for polynomial computations. https://www.singular.uni-kl.de (2025).


\bibitem{Ebe} Ebeling, W.  and Gusein-Zade S.: On the index of a vector field at an isolated singularity, The Arnoldfest (Toronto, ON, 1997), 141-152, Fields Inst. Commun.,
24, Amer. Math. Soc., Providence, RI, 1999.

\bibitem{Ebe2} Ebeling, W.  and Gusein-Zade S.: Radial index and Euler obstruction of a 1-form on a singular
variety, Geom. Dedicata 113 (2005), 231-241.





\bibitem{a08}  G\'omez-Mont, X.: An algebraic formula for the index of a vector field on a hypersurface with an isolated singularity, J. Algebraic Geom. 7, 731-752, 1998.

\bibitem{GSV} G\'omez-Mont,  X.,  Seade, J.  and  Verjovsky A.: The index of a holomorphic flow with an isolated singularity, Math. Ann. 291, 737-751, 1991.

\bibitem{ggel} Greuel, G.-M.: Dualität in der lokalen Kohomologie isolierter Singularitäten. Math. Ann. 250 (1980), 157-173.
 
\bibitem{KS} Khanedani, B.  and Suwa, T.: First variation of holomorphic forms and some applications, Hokkaido Math. J. 26 (1997), 323 - 335. 


\bibitem{King} King, H.  and Trotman, D.: Poincare-Hopf Theorems on Singular Spaces, Proceedings of the London Math. Soc., (2014) 108 (3), pp. 682-703.


\bibitem{LS01} Lehmann D. and Suwa T.: Generalization of variations and Baum-Bott residues for holomorphic foliations on singular varieties, Intern. J. Math. 10 (1999), 367-384.
 
\bibitem{a07} Lehmann, D.,  Soares, M. G.   and  Suwa, T.: On the index of a holomorphic vector field tangent to a singular variety, Bol. Soc. Bras. Mat. 26, pp. 183-199 (1995).


\bibitem{DSM} Machado, D.: Residue formula for logarithmic foliations along a divisor with isolated singularities and applications, Trans. Amer. Math. Soc., Volume 378 , 2923-2942, 2024.


\bibitem{DSM2} Machado, D.: The Schwartz index and the residue of singularities of logarithmic foliations along a hypersurface with isolated singularities. Article submitted, 2024.

\bibitem{JM} Milnor, J.: Singular points of complex hypersurfaces, Ann. of Math. Studies 61, Princeton Univ. Press, (1968).


\bibitem{Poin} Poincar\'e, H.: Sur l'integration algebrique des equations diff\'erentielles du premier ordre et du premier degr\'e, Rendiconti del Circolo Matematico di Palermo \textbf{ 5}, 161-191 (1891).

\bibitem{SSS1} Schwartz, M.-H.: Classes caractéristiques définies par une stratification d'une variété analytique complexe, C.R.Acad. Sci. Paris 260 (1965), 3262-3264,
3535-3537.

\bibitem{SSS2} Schwartz, M.-H.: Champs radiaux sur une stratification analytique complexe, Travaux en cours 39, Hermann, Paris, 1991.

 
\bibitem{SeaSuw1} Seade, J. and Suwa, T.:  A residue formula for the index of a holomorphic flow, Math. Ann. 304, 621-634 (1996).

\bibitem{SeaSuw2} Seade, J. and Suwa, T.: Residues and topological invariants of singular holomorphic
foliations, Internat. J. Math. 8 (1997), 825-847.


\bibitem{SoaAnn} Soares, M. G.: Projective varieties invariant by one-dimensional foliations, Ann. of Math. (2) 152 (2000), 369-382.

\bibitem{Soa}  Soares, M. G.:  On the geometry of Poincaré's problem for one-dimensional projective foliations. An. Acad. Brasil. Ciênc. 73 , 475-482 (2001).

\bibitem{Suw2} Suwa, T.:  Indices of vector fields and residues of singular holomorphic foliations, Actualit\'es Math\'ematiques, Hermann \'Editeurs des Sciences et des Arts (1998).

\bibitem{henrik} Vosegaard, H.: Generalized Tjurina numbers of an isolated complete intersection singularity, Math. Ann. {\bf 329} (2004), no.2, 197-224; MR2060360.

\end{thebibliography}
\end{document}